\documentclass{ifacconf}

\usepackage{graphicx}      
\usepackage{natbib}        
\usepackage{url}

\usepackage[centertags]{amsmath}	
\begin{document}

\begin{frontmatter}
\title{Modelling and Optimal Control of a Docking Maneuver with an Uncontrolled Satellite}

\author[First]{Johannes Michael} 
\author[Second]{Kurt Chudej} 
\author[Third]{J\"urgen Pannek}

\address[First]{Universit\"{a}t der Bundeswehr, Fakult\"{a}t f\"ur Luft- und Raumfahrttechnik, Institut f\"{u}r Mathematik und Rechneranwendung,
   M\"{u}nchen, Germany (e-mail: johannes.michael@unibw.de)}
\address[Second]{Universit\"{a}t Bayreuth, Fakult\"{a}t f\"{u}r Mathematik, Physik und Informatik, Lehrstuhl f\"ur Ingenieurmathematik,
   Bayreuth, Germany (e-mail: kurt.chudej@uni-bayreuth.de)}
\address[Third]{Universit\"{a}t der Bundeswehr, Fakult\"{a}t f\"ur Luft- und Raumfahrttechnik, Institut f\"{u}r Mathematik und Rechneranwendung,
   M\"{u}nchen, Germany (e-mail: juergen.pannek@unibw.de)}

\begin{abstract}                
Capturing disused satellites in orbit and their controlled reentry is the aim of the DEOS space mission. Satellites that ran out of fuel or got damaged pose a threat to working projects in orbit. Additionally, the reentry of such objects endangers the population as the place of impact cannot be controlled anymore. This paper demonstrates the modelling of a rendezvous szenario between a controlled service satellite and an uncontrolled target. The situation is modelled via first order ordinary differental equations where a stable target is considered. In order to prevent a collision of the two spacecrafts and to ensure both satellites are docked at the end of the maneuver, additional state constraints, box contraints for the control and a time dependent rendezvous condition for the final time are added. The problem is formulated as an optimal control problem with Bolza type cost functional and solved using a full discretization approach in AMPL/IpOpt. Last, simulation results for capturing a tumbling satellite are given.
\end{abstract}


\end{frontmatter}

\section{Introduction}
\label{Section:introduction}
According to NASA's \textit{Orbital Debris Program}, see \cite{NASA:FAQ}, approximately 19.000 known objects with diameters larger than ten centimeters and estimated 500.000 objects between one and ten centimeters are currently in orbit. Due to their average speed of 10 km/s even such small objects may cause serious damage to working satellites or the international space station. Additionally, with each collision the number of debris increases rapidly as the cloud of around 600 fractions of the `Cosmos 2251' and `Iridium' crash in 2009, see \cite{GUARDIAN}. \\
Aside from orbital collisions, an uncontrolled reentry into the atmosphere poses a severe threat. As the entrance corridor is highly difficult to precalculate, with every crash pieces of beryllium or titan may reach the earth's surface and hit populated areas. In October 2011 the german satellite `Rosat' crashed uncontrolled and even a few hours prior to impact experts were unsure where fractions will strike the earth, see \cite{ROSAT}. \\
The planned german technology demonstration mission DEOS (\textbf{De}ut\-sche \textbf{O}r\-bi\-ta\-le \textbf{S}ervicing Mission) takes the next step towards capturing and controlling disused satellites. Considered as a feasibility study, DEOS consists of two spacecrafts, one acting as uncontrollable target and one as servicer. The aims of the mission are to demonstrate docking and berthing with a noncooperating target and the controlled reentry into atmosphere in coupled condition, see \cite{DEOS:2010}. \\
This paper presents a model to calculate the optimal control schemes for a DEOS capturing mission. In particular, the equations for a rendezvous maneuver to an uncontrollable target are derived and solved for one example with a direct approach. The presented results show the feasibility of such a mission. Additionally, reference trajectories and control schemes for different situations are obtained which are extendable to the case of a tumbling target.\\
The following section deals with the derivation of the optimal control problem. The equations of motion for the relative position and orientation, the terminal condition, the state and control contraints and the cost functional are presented. Afterwards a direct approach is applied to solve the optimal control problem. In the fourth part the numerical results for a flyaround maneuver with a stable spinning target are presented and discussed. Thereafter a short introduction in simulating the trajectories in a virtual reality is given. An outlook of possible extensions concludes the paper.

\section{Setup of the rendezvous model}
\label{Section:setup}
The motion of each spacecraft consists of two independent 3D subsystems, the relative position and the relative orientation. The equations of motion for the relative position are the so called Clohessy-Wilshire equations. Further information on orbit analysis and design can be found in \cite{SCHI:11}. Note that throughout this work we shorten notation by omitting the time dependency if no specific time instant is considered.

\subsection{Relative position}
Assuming circular orbits the subsystem of the relative position is derived via Kepler's equation for the two body problem and a linear taylor expansion, see, e.g., \cite{kyle}. As a result, a second order ordinary differential equation system for the position of a body relative to a moving reference point in orbit is obtained. Throughout this work, we consider the uncontrolled target to be that reference point. Apart from the standard state variables $x$, $y$ and $z$ denoting the relative distance in radial, direction of movement and out of orbit component respectively, we introduce the additional state variables $v_x$, $v_y$ and $v_z$ denoting the respective velocities to transform the system into six first order differential equations. 
  \begin{align}
  	 \label{eq:cw_eq1} 
    \dot{x} &= v_x \displaybreak[0] \\ 
    \dot{y} &= v_y \displaybreak[0] \\
    \dot{z} &= v_z \displaybreak[0] \\
    \dot{v}_x &= 2n v_y + 3n^{2} x  + \frac{u_{x}}{m} \displaybreak[0] \\
    \dot{v}_y &= - 2n v_x + \frac{u_{y}}{m} \displaybreak[0] \\
    \dot{v}_z &= - n^{2}z + \frac{u_{z}}{m} \label{eq:cw_eq6}
  \end{align}
Here $m$ denotes the mass of the servicer spaecraft and $n$ is the mean motion of the reference point, i.e. the angular velocity in orbit. For the dynamics it can be seen that the $z$ coordinate is decoupled from the $x$-$y$-system, i.e. the out of orbit motion is independent of the motion in the orbit plane. If no controls $u_x, \;u_y, \;u_z$ are applied an analytical solution is known, again see \cite{kyle}. Yet for non vanishing controls no analytical solution is known. For this reason, we evaluate the dynamics numerically.

\subsection{Orientation}
In space mission modelling one commonly uses the quaternion representation to derive the equations for the orientation of a satellite. Quaternions are a special parametrization of the Euler--axis/angle description which is nonsingular and continuous. A quaternion can be considered as a four dimensional vector $\boldsymbol{q} := [q_1,\,q_2,\,q_3,\,q_4]^\top$. The first three components are called the vector part and the fourth is the scalar component. The denomination results from the fact that a quaternion represents a rotation around an axis in three dimensional space, i.e. the vector part, by a specific angle given by the scalar component. To get a unique description of the orientation, the quaternions are supposed to be of length one, so called unit quaternions. Further information on quaternions and their relation to other rotation types can be found in \cite{tewari}.\\
The first derivative of quaternions is given by
\begin{align*} 
  \dot{\boldsymbol{q}} = \frac{1}{2}
  \begin{bmatrix}
  \boldsymbol{\omega} \\
  0
  \end{bmatrix}
  \otimes \boldsymbol{q} ,
\end{align*}
where $\boldsymbol{\omega} := [\omega_x, \, \omega_y,\, \omega_z]^\top$ is the vector of the angular velocities and $\otimes$ represents the quaternion multiplication, see also \cite{stevens:03} for a derivation. In its components the dynamics of the quaternions are of the form
  \begin{align} \label{eq:quat_der_expl}
    \begin{bmatrix}
      \dot{q}^{\alpha}_{1}\\
      \dot{q}^{\alpha}_{2}\\
      \dot{q}^{\alpha}_{3}\\
      \dot{q}^{\alpha}_{4}
    \end{bmatrix} &= \frac{1}{2}
    \begin{bmatrix}
      0 & \omega^{\alpha}_{z} & -\omega^{\alpha}_{y} & \omega^{\alpha}_{x} \\
      -\omega^{\alpha}_{z} & 0 & \omega^{\alpha}_{x} & \omega^{\alpha}_{y} \\
      \omega^{\alpha}_{y} & -\omega^{\alpha}_{x} & 0 & \omega^{\alpha}_{z} \\
      -\omega^{\alpha}_{x} & -\omega^{\alpha}_{y} & -\omega^{\alpha}_{z} & 0
    \end{bmatrix}
    \begin{bmatrix}
      q^{\alpha}_{1}\\
      q^{\alpha}_{2}\\
      q^{\alpha}_{3}\\
      q^{\alpha}_{4}
    \end{bmatrix} \;
    \alpha \in \{S, T\},
  \end{align}\\
where $S$ and $T$ denote the service and target spacecraft respectively. \\
Note that the angular velocity of a tumbling object is typically not constant with respect to time. Denoting the inertia tensor with $\boldsymbol{J}$ we incorporate that aspect using Euler's gyroscopic equation
\begin{align*}
  \dot{\boldsymbol{J}} \cdot \boldsymbol{\omega} + \boldsymbol{J} \cdot \dot{\boldsymbol{\omega}} + \boldsymbol{\omega} \times \left(\boldsymbol{J} \cdot \boldsymbol{\omega}\right) = \boldsymbol{m},
\end{align*}
which provides the needed information, see, e.g., \cite{chobotov}. Under the assumptions that mass distribution is constant over time, i.e. $\dot{\boldsymbol{J}}=0$, and the regarded system coincides with the principal axis of the body, that is $J_{ik} = 0$ for all $i \not = k$, we obtain the set of first order equations 
  \begin{align}
    \dot{\omega}^{S}_{x} &= \frac{1}{J^{S}_{xx}}\left(\omega^{S}_{y}\omega^{S}_{z}\left(J^{S}_{yy} - J^{S}_{zz}\right) + m_x\right) \label{eq:gyro1} \displaybreak[0] \\
    \dot{\omega}^{S}_{y} &= \frac{1}{J^{S}_{yy}}\left(\omega^{S}_{x}\omega^{S}_{z}\left(J^{S}_{zz} - J^{S}_{xx}\right) + m_y\right) \displaybreak[0] \\
    \dot{\omega}^{S}_{z} &= \frac{1}{J^{S}_{zz}}\left(\omega^{S}_{x}\omega^{S}_{y}\left(J^{S}_{xx} - J^{S}_{yy}\right) + m_z\right) \displaybreak[0] \\
    \dot{\omega}^{T}_{x} &= \frac{1}{J^{T}_{xx}}\left(\omega^{T}_{y}\omega^{T}_{z}\left(J^{T}_{yy} - J^{T}_{zz}\right)\right) \displaybreak[0] \\
    \dot{\omega}^{T}_{y} &= \frac{1}{J^{T}_{yy}}\left(\omega^{T}_{x}\omega^{T}_{z}\left(J^{T}_{zz} - J^{T}_{xx}\right)\right) \displaybreak[0] \\
    \dot{\omega}^{T}_{z} &= \frac{1}{J^{T}_{zz}}\left(\omega^{T}_{x}\omega^{T}_{y}\left(J^{T}_{xx} - J^{T}_{yy}\right)\right). \label{eq:gyro6}
  \end{align}
where $\boldsymbol{m} := [m_x,\;m_y,\;m_z]^\top$ defines the momentum vector. Note that no momentum $\boldsymbol{m}$ is applied on the right hand side of the target satellite as it is supposed to be uncontrolled.\\
With these equations of motion we obtain a system of twenty first order ordinary differential equations \eqref{eq:cw_eq1}--\eqref{eq:gyro6}. Here, we denote the combined state vector by
\begin{align*}
  \boldsymbol{x} := [&x, \;\, y, \;\, z, \;\, v_x, \;\, v_y, \;\, v_z, \;\, \omega^{S}_{x}, \;\, \omega^{S}_{y}, \;\, \omega^{S}_{z}, \;\, \omega^{T}_{x}, \;\, \omega^{T}_{y}, \nonumber\\
  &  \;\, \omega^{T}_{z}, \;\, q^{S}_{1}, \;\, q^{S}_{2}, \;\, q^{S}_{3}, \;\, q^{S}_{4}, \;\, q^{T}_{1}, \;\, q^{T}_{2}, \;\, q^{T}_{3}, \;\,  q^{T}_{4} ]^\top 
\end{align*}
and the control vector by
\begin{align*}
	\hat{\boldsymbol{u}} := \left[u_x, \;\, u_y, \;\, u_z, \;\, m_x, \;\, m_y, \;\, m_z \right]^\top.
\end{align*}

\subsection{Terminal condition}
To model the docking of the two spacecraft some terminal constraints have to be introduced. We define a docking point for each satellite in their body fixed coordinate system, $\boldsymbol{d}^{S} := [d^{S}_{1},\;d^{S}_{2},\;d^{S}_{3}]^\top$ and $\boldsymbol{d}^T:= [d^{T}_{1},\;d^{T}_{2},\;d^{T}_{3}]^\top$, whose local positions and velocities have to match at the end of the maneuver. Similar conditions can be found in, e.g., \cite{kyle} and \cite{boyarko}. To get a common representation, all components are converted to the earth centered inertial coordinate system. We denote the corresponding transformation matrix from the rotated body fixed coordinates to the inertial system by $R^\top$ and add the index $\mathcal{E}$ to variables defined in the earth central inertial coordinate system. Computing the transformation matrix from the components of the related quaternion, see, e.g., \cite{tewari}, we obtain
\begin{small}
\begin{align} \label{eq_rotmat}
  R = 
  \left[
    \begin{array}{ccc}
	q^{2}_{1} - q^{2}_{2} - q^{2}_{3} + q^{2}_{4} & 2\left(q_{1}q_{2} + q_{3}q_{4}\right) & 2\left(q_{1}q_{3} - q_{2}q_{4}\right) \\
	2\left(q_{1}q_{2} - q_{3}q_{4}\right) & -q^{2}_{1} + q^{2}_{2} - q^{2}_{3} + q^{2}_{4} & 2\left(q_{2}q_{3} + q_{1}q_{4}\right)\\
	2\left(q_{1}q_{3} + q_{2}q_{4}\right) & 2\left(q_{2}q_{3} - q_{1}q_{4}\right) & -q^{2}_{1} - q^{2}_{2} + q^{2}_{3} + q^{2}_{4}
    \end{array}
    \right]
\end{align}
\end{small} 
Multiplication with this matrix performs the transformation of a vector from the unrotated coordinate system to the rotated one. The inverse transformation, as we need it in this case, can be done by multiplication with the transposed matrix of R. Applying the transformation for both the service and the target satellites, the local difference of the docking point is given by
\begin{align} \label{rho_komplett}
  \boldsymbol{d}^{S}_{\mathcal{E}} - \boldsymbol{d}^{T}_{\mathcal{E}} = \boldsymbol{r} + {R^{S}}^\top \boldsymbol{d}^{S} - {R^{T}}^\top \boldsymbol{d}^{T},
\end{align}
with $\boldsymbol{r} := [x,\;y,\;z]^\top$ being the cartesian distance of the two spacecrafts. Similarly, the angular velocity of each satellite in the earth centered inertial frame is given by
\begin{align*}
  \boldsymbol{\omega}^{\alpha}_\mathcal{E} = {R^\alpha}^\top \boldsymbol{\omega}^{\alpha} - 
  \begin{bmatrix}
  0 \\
  0 \\
  n
  \end{bmatrix} \qquad \alpha \in \left\{S,T\right\}.
\end{align*}	
Hence, we obtain the difference of the velocity via
\begin{align} \label{rhopunkt_komplett}
  \dot{\boldsymbol{d}}^{S}_{\mathcal{E}} - \dot{\boldsymbol{d}}^{T}_{\mathcal{E}} = \dot{\boldsymbol{r}} + \boldsymbol{\omega}^{S}_{\mathcal{E}} \times \left( {R^{S}}^\top \boldsymbol{d}^{S}\right) - \boldsymbol{\omega}^{T}_{\mathcal{E}} \times \left( {R^{T}}^\top \boldsymbol{d}^{T}\right).
\end{align}
To simplify this equation we assume that the angular velocity and the orientation of both spacecrafts coincide at the end of the time interval, i.e.,
  \begin{align} \label{eq:ass_final_rot}
     \boldsymbol{q}^{T}(t_f) - \boldsymbol{q}^{S}(t_f) &= \boldsymbol{0} \\
    \boldsymbol{\omega}^{T}(t_f) - \boldsymbol{\omega}^{S}(t_f) &= \boldsymbol{0}.\label{eq:ass_final_angvel}
  \end{align}
Note that the docking point definition is required to be conform with these assumptions. To illustrate this point, consider the docking points to be located at the back of the satellites. Then,  since the satellites would have to overlap, there exists no physically possible configuration where these points and the satellite orientations coincide.\\
Assuming \eqref{eq:ass_final_rot} and \eqref{eq:ass_final_angvel} to hold we obtain
\begin{align*}
 R^{T}(t_f) = R^{S}(t_f) = R(t_f) \\
 \omega^{T}_{\mathcal{E}}(t_f) = \omega^{S}_{\mathcal{E}}(t_f) = \omega_{\mathcal{E}}(t_f)
\end{align*}
and the docking conditions \eqref{rho_komplett} and \eqref{rhopunkt_komplett} simplify to
\begin{align*}
    R(t_f)^\top \left( \boldsymbol{d}^{T} - \boldsymbol{d}^{S} \right) - \boldsymbol{r} &= \boldsymbol{0} \\
    \boldsymbol{\omega}_\mathcal{E}(t_f) \times R(t_f)^\top (\boldsymbol{d}^{T} - \boldsymbol{d}^{S}) - \dot{\boldsymbol{r}} &= \boldsymbol{0}.
\end{align*}

\subsection{State and control constraints}
\label{Subsection:constraints}
As pointed out in the previous subsection, the case of colliding satellites has to be treated. To exclude such an occurance, we first introduce a spherical safety area around each satellite with radius $r^{S} > 0$ and $r^{T} > 0$ respectively. In our case, the usage of spheres is reasonable as the geometry of the satellites is supposed to be simple. Based on the safety areas, we add the state constraint 
\begin{align*}
x^2 + y^2 + z^2 \geq (r^{S} + r^{T})^2
\end{align*}
to the optimal control problem which resolves the possible collision issue. Note that the safety area needs to be chosen carefully since the state constraint together with the terminal condition may exclude the existence of a solution. \\
Apart from the states we also consider the controls to be constrained. In particular, we suppose that the thrust vector $\boldsymbol{u} := [u_x,\; u_y, \; u_z]^\top$ is bounded via
\begin{align*}
	u^2_x + u^2_y + u^2_z \leq u_{\max}
\end{align*}
and each momentum is bounded via
\begin{align*}
	\left| m_i \right| \leq m_{\max} \quad i=x, y, z.
\end{align*}
The choice of the limitations in the momentum control is motivated by the fact that the satellites attitude control system can be realized via control moment gyroscopes. These can cause a maximal momentum $m_{\max}$ around each axis. The bound of the thrust vector is due to future extensions of the model where only one thruster is installed and the spacecraft has to rotate to accelerate in a specific direction. \\
As the thrusters are mounted on the satellite their direction depends on the current attitude of the spacecraft. Therefore we transform the calculated thrust $u_x, \; u_y$ and $u_z$ from earth centered coordinates into the body fixed system. The transformation is done by multiplication with the rotation matrix $R$ defined in \eqref{eq_rotmat} and reveals
\begin{align*}
 \begin{bmatrix}
  u_1 \\
  u_2 \\
  u_3
 \end{bmatrix} := R
 \begin{bmatrix}
  u_x \\
  u_y \\
  u_z
 \end{bmatrix}
\end{align*}
where $u_1, \; u_2$ and $u_3$ denote the controls in the body fixed coordinates. Note that as $R$ is orthogonal the length of $\boldsymbol{u}$ remains unchanged by this transformation, and hence the upper bound $u_{\max}$ on the controls in the body fixed and the earth centered coordinates are identical.

\subsection{Cost functional}
Throughout this work we consider {\it minimizing} the following Bolza type cost functional over the set of allowed control functions $\hat{\boldsymbol{u}}(t) := [\boldsymbol{u}(t), \; \boldsymbol{m}(t)]^\top$ and terminal times $t_f$. We define the cost functional as follows
\begin{align} \label{ZF}
  	J(\hat{\boldsymbol{u}}(t),t_f) = l_{t_f} t_f + \int^{t_f}_{0}l_{u}\left\| \boldsymbol{u}(t) \right\|_2^2 + l_{m} \left\| \boldsymbol{m}(t) \right\|_2^2dt.
\end{align}
This functional represents a combination of the terminal time, the position controls and the orientation controls. The non-negative constants $l_{t_f}$, $l_u$ and $l_m$ can be used to force scenarios such as time-optimalty or minimal control effort. Note that the $2$--norm $\| \cdot \|_2$ has to be used if results shall hold for both controls $(u_1, u_2, u_3)^\top$ and $(u_x, u_y, u_z)^\top$.

\section{Discretization and implementation}
\label{Section:discretization}
To solve the optimal control problem
\begin{align}
	&\mbox{Minimize } \quad J(\hat{\boldsymbol{u}}(t),t_f) \nonumber\\ 
	 &\text{subject to the dynamics} \quad \eqref{eq:cw_eq1} - \eqref{eq:gyro6} \nonumber\\
	 &\text{with initial and terminal conditions} \nonumber\\
	  & \qquad \boldsymbol{x}(0) = \boldsymbol{x}_0 \displaybreak[0] \nonumber\\
	  & \qquad \boldsymbol{0} = \boldsymbol{q}^{T}(t_f) - \boldsymbol{q}^{S}(t_f) \displaybreak[0] \nonumber\\
	  & \qquad \boldsymbol{0} = \boldsymbol{\omega}^{T}(t_f) - \boldsymbol{\omega}^{S}(t_f) \displaybreak[0] \nonumber\\
	  & \qquad \boldsymbol{0} = R(t_f)^\top \left( \boldsymbol{d}^{T} - \boldsymbol{d}^{S} \right) - \boldsymbol{x}  \displaybreak[0] \nonumber\\
	  & \qquad \boldsymbol{0} = \boldsymbol{\omega}_\mathcal{E}(t_f) \times R(t_f)^\top (\boldsymbol{d}^{T} - \boldsymbol{d}^{S}) - \dot{\boldsymbol{x}} \displaybreak[0] \label{eq:const1}\\
	    & \text{and constraints} \nonumber\\
	  & \qquad u_{\max} \geq u^2_x(t) + u^2_y(t) + u^2_z(t) \quad \forall t \in [0, t_f] \displaybreak[0] \nonumber \\
	  & \qquad m_{\max} \geq \left| m_i(t) \right| \quad i=x, y, z \quad \forall t \in [0, t_f] \displaybreak[0] \nonumber \\
	  & \qquad x^2(t) + y^2(t) + z^2(t) \geq (r^{S} + r^{T})^2 \quad \forall t \in [0, t_f] \label{eq:const2}
\end{align}
a direct approach can be applied. Using a full discretization with fixed step size $\Delta t = t_f/N$ where $N$ is the number of steps, the discretized dynamics of the state $\boldsymbol{x}$ and the discretized constraints \eqref{eq:const1}, \eqref{eq:const2} form the constraints of the resulting finite nonlinear program. The optimization variable now contains not only the discretization of the control $\hat{\boldsymbol{u}}$ and the free terminal time $t_f$, but also the discretization of the entire state trajectory. With the exponent $k$ we denote the state of a variable at the $k$th discretization point, i.e. $x^k=\boldsymbol{x}(k\cdot \Delta t)$.\\
Since the gyroscopic equation is stiff the implicit trapezoidal method
\begin{align*}
	\frac{\boldsymbol{x}^{k+1} - \boldsymbol{x}^k}{\Delta t} = \frac{1}{2}(\boldsymbol{f}\left(\boldsymbol{x}^k, \boldsymbol{u}^k,t^k\right) + \boldsymbol{f}\left(\boldsymbol{x}^{k+1}, \boldsymbol{u}^{k+1},t^{k+1}\right)) \\
	k = 0,1,...,N-1
\end{align*}
is used to discretize the dynamics of the state. Note that the components of the angular velocity would soar even for very small step-sizes if, e.g., Euler's method was used. We implemented the discretized optimal control problem in the modelling language AMPL and solved it using the interior point optimizer IpOpt. For details on AMPL and IpOpt see, e.g., \cite{ampl} and \cite{WB2006}. For the discretized problem, the constraints are defined pointwise at the discretization points. To avoid constraint violations at intermediate time instances, adaptation techniques for the discretization grid may be used, see, e.g., \cite{BH1998}, or the bounds may be tightened.

\section{Numerical Results}
\label{Section:numerical results}

To illustrate our results we consider a flyaround maneuver, that is a situation where the service satellite is located on the wrong side of the target pointing with its docking interface towards it. Here, we assume that the docking interface of the servicer is located on the 'front' whereas the docking point of the target is at its 'back', e.g. the targets thruster where some kind of hook clasps. The radius of the safety area for each satellite is chosen to be of size one. For both satellites we assume that they are axially symmetric with respect to the $y$ axis which is justified for the regarded satellite geometries. Note that the choice of the docking points enables us to use the simplified equations for the terminal condition \eqref{eq:ass_final_rot} and \eqref{eq:ass_final_angvel} and creates no impossible final states.

Upon start of the maneuver, we suppose that the servicer is in a non-rotated initial state and its angular velocity is zero. The target also starts in a non-rotating initial state but in contrast to the servicer it is rotating with a constant angular velocity of $3^\circ/s$ around its $y$ axis resulting in a stable motion with respect to time. The initial values are as follows:
\begin{align*}
\begin{bmatrix}
x_0 \\
y_0 \\
z_0 
\end{bmatrix} &= 
\begin{bmatrix}
0 \\
3 \\
0 
\end{bmatrix}, \qquad
\begin{bmatrix}
v_{x0} \\
v_{y0} \\
v_{z0} \\
\end{bmatrix} = 
\begin{bmatrix}
0 \\
0 \\
0 
\end{bmatrix}, \displaybreak[0] \\
\boldsymbol{q}^{S} &=
\begin{bmatrix}
0 \\
0 \\
1 \\
0 
\end{bmatrix}, \qquad \quad
\boldsymbol{\omega}^{S} =
\begin{bmatrix}
0 \\
0 \\
0 
\end{bmatrix}, \displaybreak[0]  \\
\boldsymbol{q}^\top &=
\begin{bmatrix}
0 \\
0 \\
0 \\
0 
\end{bmatrix}, \qquad \quad
\boldsymbol{\omega}^\top =
\begin{bmatrix}
0 \cdot 0.017453 \\
3 \cdot 0.017453 \\
0 \cdot 0.017453 
\end{bmatrix}.
\end{align*}
Geometrically, the above conditions force the servicer to fulfill three tasks: It has to fly around the target, turn around to align the docking points and adopt the rotation such that the angular velocities coincide. 

Within the cost functional \eqref{ZF} we consider the weights $l_{t_f} = 1$, $l_u = 1$ and $l_m = 1$, i.e., the terminal time and the control costs are weighted equally. The remaining constants of the optimal control problem are displayed in Table~\ref{tb:constants}. 
\begin{table}[!htb]
\begin{center}
\begin{tabular}{lcl}
variable & value & description \\\hline
$a$ & 7071000 & orbit radius [m]\\
$GM$ & $398 \cdot 10^{12}$ & gravitational constant [$N (m/kg)^2$]\\ 
$n$ & $\sqrt{\frac{GM}{a^3}}$ & mean motion [1/s] \\
$m$ & 200 & satellite mass [kg] \\
$u_{\max}$ & 0.15 & maximum thrust [N]\\
$m_{\max}$ & 1 & maximum torgue [Nm] \\
$J^T_{xx}$ & 1000 & Targ. angular mass around x [$kg/m^2$]\\
$J^T_{yy}$ & 2000 & Targ. angular mass around y [$kg/m^2$]\\
$J^T_{zz}$ & 1000 & Targ. angular mass around z [$kg/m^2$]\\
$J^{S}_{xx}$ & 2000 & Serv. angular mass around x [$kg/m^2$]\\
$J^{S}_{yy}$ & 5000 & Serv. angular mass around y [$kg/m^2$]\\
$J^{S}_{zz}$ & 2000 & Serv. angular mass around z [$kg/m^2$]\\
$r^{T}$ & 1 & safety-area around target [m]\\
$r^{S}$ & 1 & safety-area around servicer [m]\\
$d^{T}$ & $\left[0 \;\; -1 \;\; 0\right]^\top$ & docking point target [m] \\
$d^{S}$ & $\left[0 \;\; 1 \;\; 0\right]^\top$ & docking point servicer [m]
\end{tabular}
\caption{Constants of optimal control problem}\label{tb:constants}
\end{center}
\end{table}

\begin{figure}[htb]
  \begin{center}
  \includegraphics[width=\linewidth]{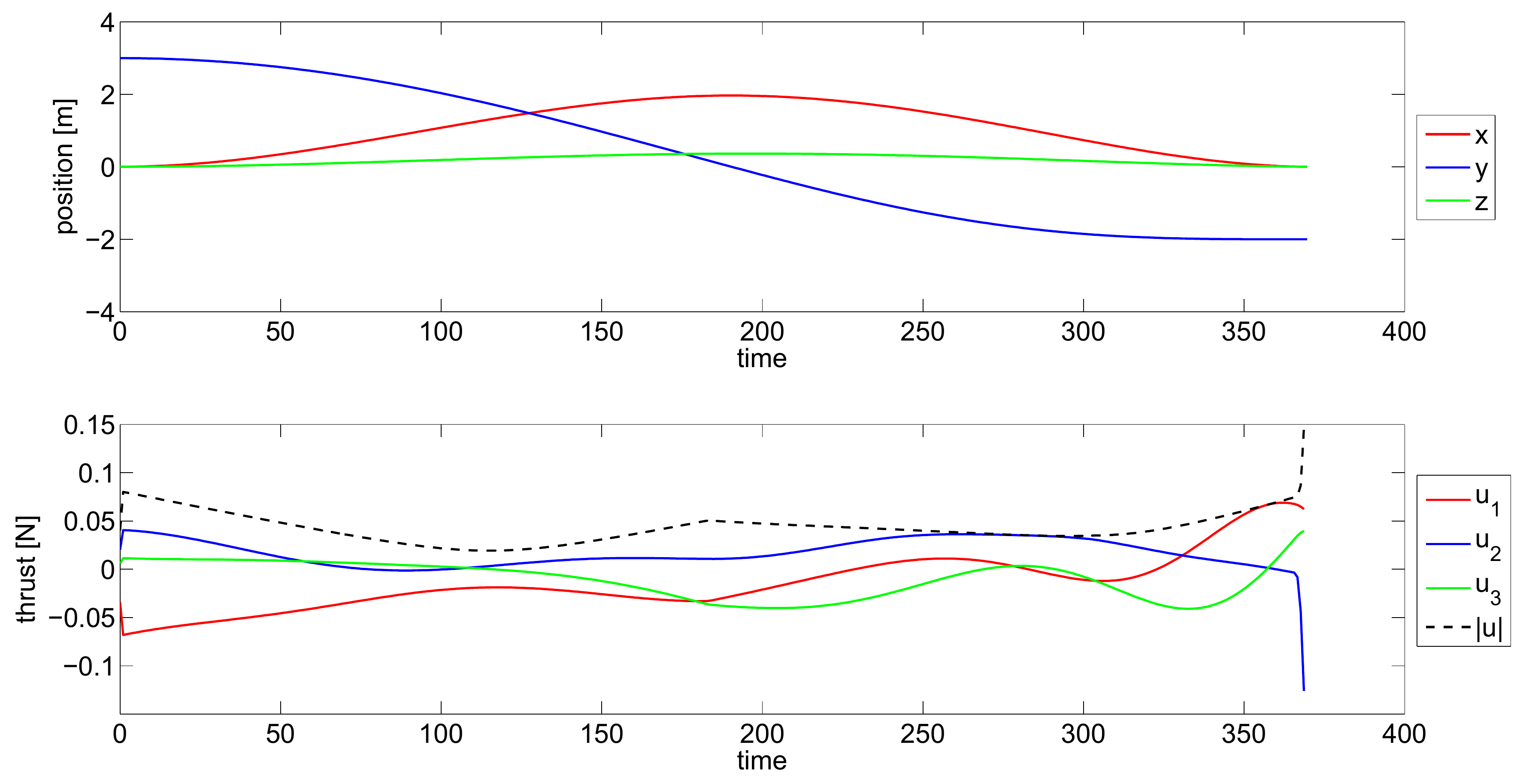}    
  \caption{Position and thruster control of the servicer} 
  \label{fig:position}
  \end{center}
\end{figure}

\begin{figure}[htb]
  \begin{center}
  \includegraphics[width=\linewidth]{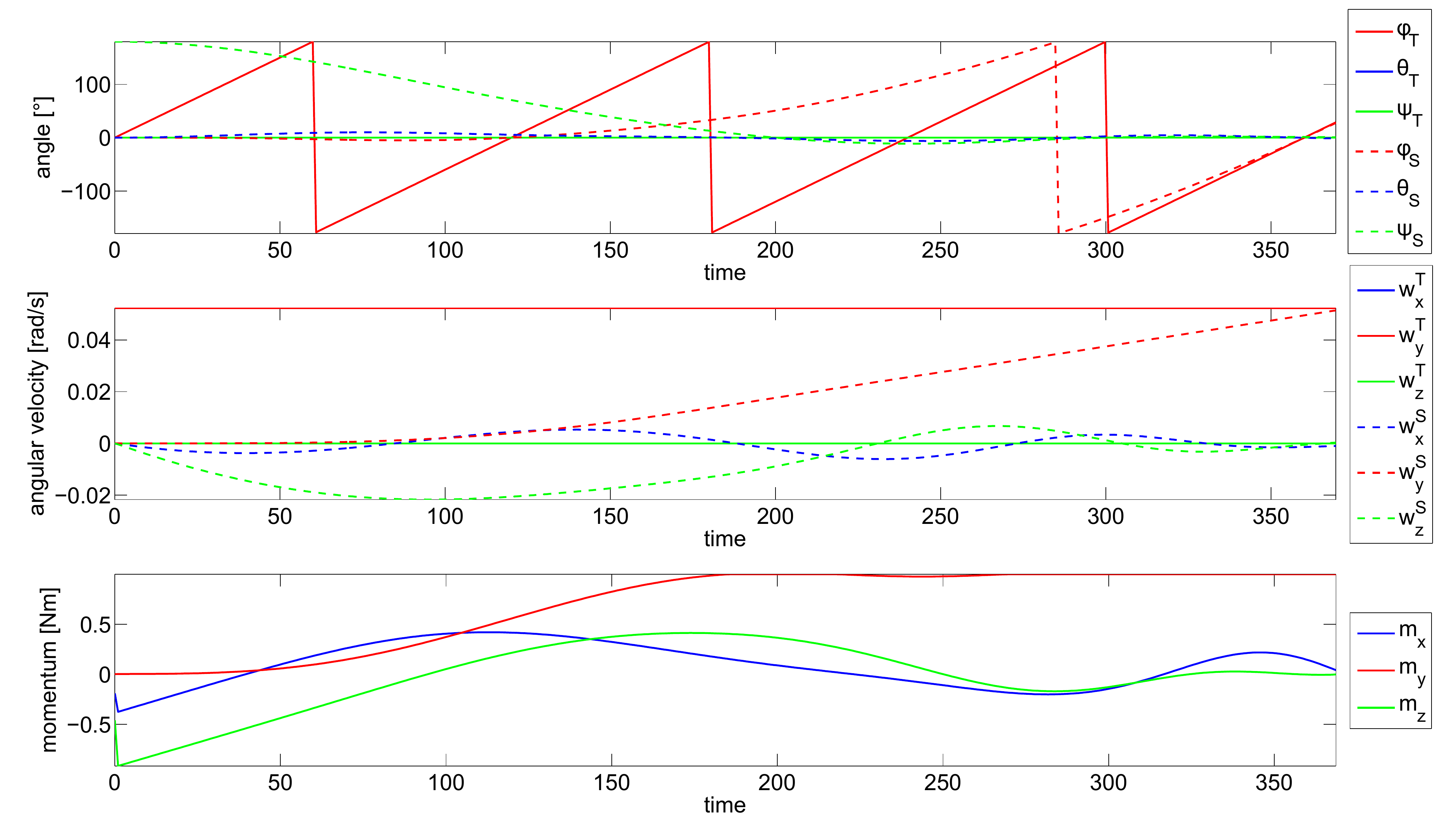}    
  \caption{Orientation, angular velocity and momentum controls of target and servicer (dashed)} 
  \label{fig:rotation}
  \end{center}
\end{figure}

To convert the optimal control problem into a finite nonlinear program we used a discretization with $N=370$ steps. The results of our calculations are displayed in Figures \ref{fig:position} and \ref{fig:rotation} which show the states and controls over time. For the optimal terminal time we obtain $t_f = 369.61$ seconds, which is about $6.16$ minutes. \\

From Figure \ref{fig:position} one observes that the position control never reaches the maximum. The fact that the allowed thrust is never applied is due to the chosen weights in the cost functional. The controls are smooth, except one point at about $190$ seconds, but there occurs no bang-bang behaviour. Because of the state constraints the satellite cannot fly to its terminal position in a straight line. In the beginning it accelerates in positive $x$ direction -- which is due to the initial orientation -- performed by a negative control of $u_1$. Note that the applied control has to be interpreted with respect to the current attitude of the satellite, cf. Section \ref{Subsection:constraints}. Upon termination of the maneuver the negative control in $y$ direction represents the breaking to prevent collision with the target and to obtain zero relative velocity. As expected, the position in $x$ and $y$ nearly produce a symmetrical arc with respect to half the maneuver time.

The orientations in Figure \ref{fig:rotation} are plotted in Euler-$yxz$-angles with the angle notation $\phi$, $\theta$ and $\psi$. One observes that the spin of the target around the $y$-axis is periodical. In the first part of the maneuver the main momentum control of the servicer is applied to perform the turning. To this end, the almost maximal momentum is applied around the $z$-axis and some additional correction in $x$ direction. This control first accelerates the angular velocity of the servicer around the turning axis and then slows it down such that upon termination of the maneuver the angular velocity around the $x$- and $z$-axis are zero. After about 50 seconds the momentum around the $y$-axis is increasing. It reaches the maximal possible amount and nearly keeps this value until the end of the maneuver. This is the control for adapting the angular velocity around $y$. The acceleration of $\omega_y$ is timed such that upon termination both the angular velocity and the orientation of servicer and target are equal.

Computing the total control costs of the thrust and the momentum control, we obtain
\begin{align*}
	u_{\text{total}} & = \Delta t \sum^{N-1}_{k=0} (u^k_x)^2 + (u^k_y)^2 + (u^k_z)^2 = 15.7662  \displaybreak[0] \\
	m_{\text{total}} & = \Delta t \sum^{N-1}_{k=0} (m^k_x)^2 + (m^k_y)^2 + (m^k_z)^2 = 295.5767.
\end{align*}
Hence, the costs of the momentum control outweigh the costs of the position control for this particular example. The optimal value of the cost functional sums up to
\begin{align*}
J = t_f + u_{\text{total}} + m_{\text{total}} = 680.9548.
\end{align*}

\begin{figure}[htb]
  \begin{center}
  \includegraphics[width=\linewidth]{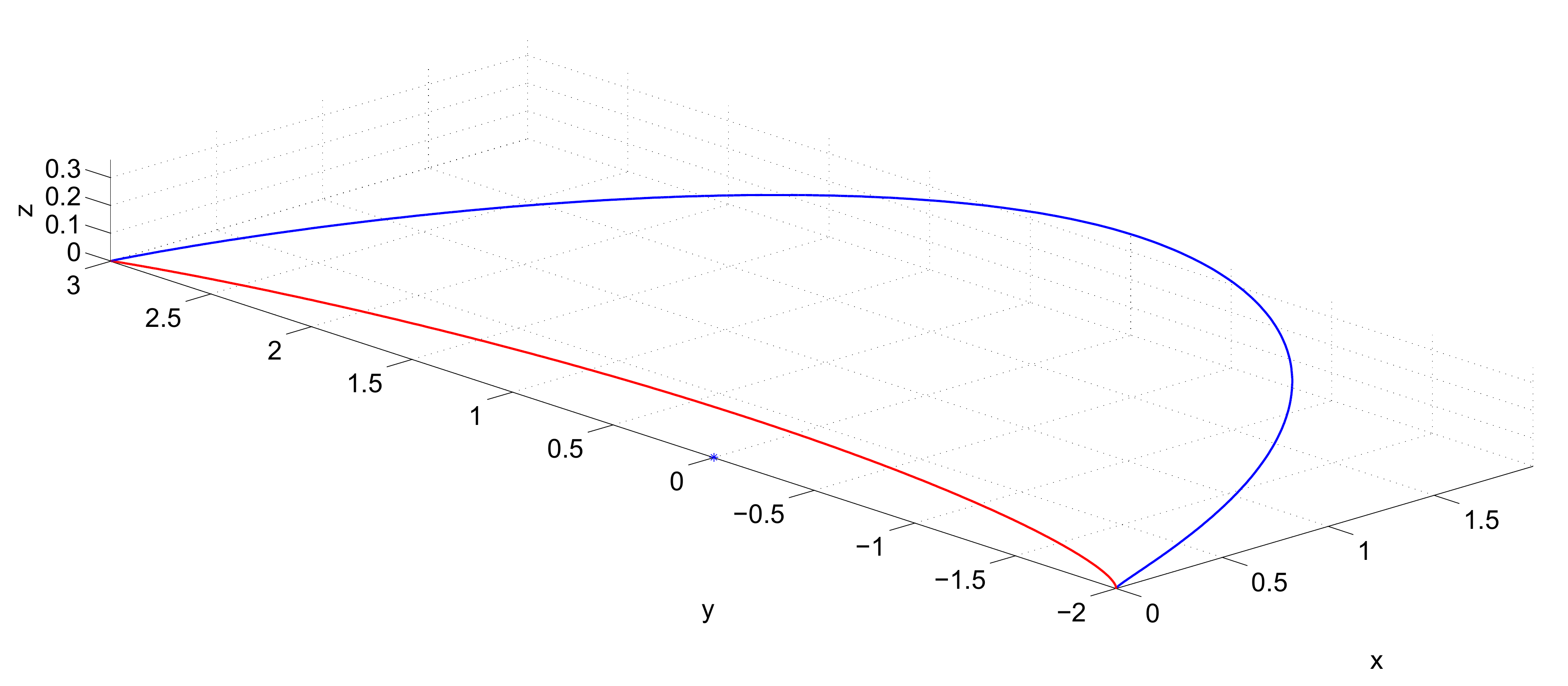}    
  \caption{Trajectory of the servicer with state constraint (blue) and without state constraint (red)} 
  \label{fig:pos3d}
  \end{center}
\end{figure}

To demonstrate that the calculated control depends on the state constraints we calculated the optimal trajectory with the same initial values omitting this constraint. In Figure \ref{fig:pos3d} one can see that with the safety areas the servicer flies a nearly circular arc around the target. Without the constraint the flight path is much more direct and the servicer would crash into the target. 

\section{Simulation}
\label{Section:simulation}
\begin{figure}[htb]
\begin{center}
\includegraphics[width=\linewidth]{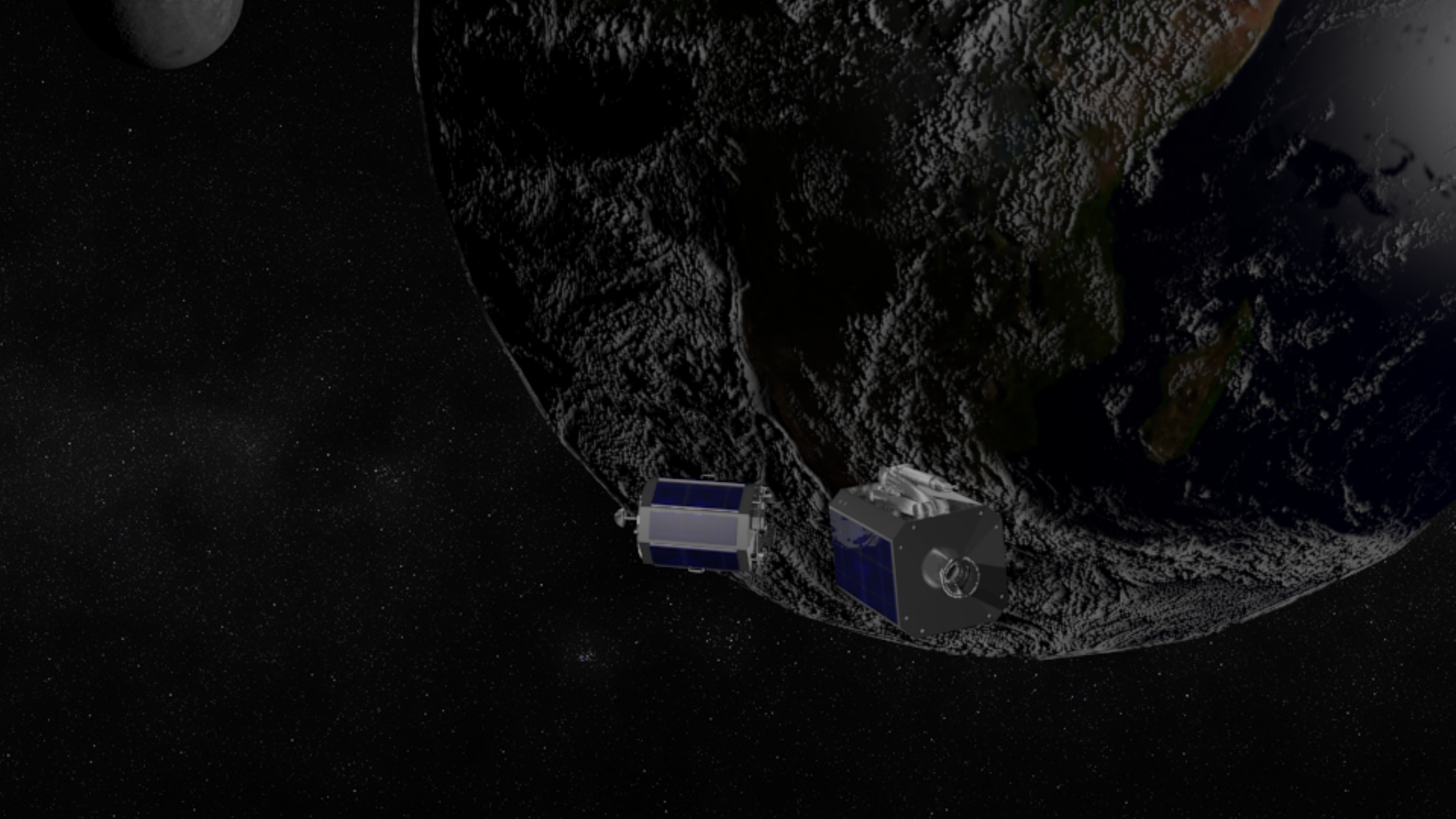}    
\caption{Rendered picture of the satellites} 
\label{fig:rendered}
\end{center}
\end{figure}

Since the computed results are hard to imagine from the plots, we developed a virtual reality where the calculated data can be imported and rendered images of the maneuver can be created. These images are then combined to a video showing the entire rendezvous maneuver. The virtual reality was created with the open source software Blender 3D. Using an additional import tool the trajectory data is bound to the corresponding satellite. Moreover, the calculated quaternions are used directly to rotate the models. For information on Blender see, e.g., \cite{blender}. The satellite models were designed based on the first design studies of the german aerospace center, see \cite{DEOS:2010}. In Figure \ref{fig:rendered} one can see an image of the two satellites in the created virtual reality.

\section{Outlook}
\label{Section:outlook}
The shown example only covers a stable rotating target satellite, i.e. all angular velocity components are constant with respect to time. Yet, the presented model can also be used to calculate optimal solutions for a tumbling target which may occur if a satellite runs out of fuel or energy to keep its orientation. This case is under investigation at the moment to handle the additional difficulty of a moving target point. Note that it is possible to calculate the maximal tumbling motion which still allows a docking maneuver to be performed. If this maximum is exceeded, the service satellite can be fitted with an additional manipulator arm to grasp the target to establish a connection.  To optimize such a so called berthing maneuver the model derived in this work has to be modified by adding a model of a robot arm, see also \cite{fehse2003automated} for information on the differences between docking and berthing. For one, future research will concern computing optimal berthing maneuvers, but also to derive low cost feedback controllers for both docking and berthing situations.

\begin{ack}
We thank Prof. Dr. K. Schilling, University of W\"urzburg, for suggesting this work and for fruitful discussions. 
\end{ack}


\end{document}